\documentclass[11pt,final]{article}

\usepackage{lineno,hyperref, xcolor}
\usepackage[margin=1in]{geometry}
\usepackage{subfig}
\usepackage{pgfplots}
\usepackage{ragged2e}

\usepackage{natbib}
 \bibpunct[, ]{(}{)}{,}{a}{}{,}%
\begin{document}

\title{Constrained Local Search for Last-Mile Routing}

\author
{
  William Cook \\ University of Waterloo\footnote{Department of Combinatorics and Optimization}\\
  bico@uwaterloo.ca
\and
Stephan Held \\ University of Bonn\footnote{Research Institute for Discrete Mathematics and Hausdorff Center for Mathematics}\\
held@dm.uni-bonn.de
\and
Keld Helsgaun\\ Roskilde University\footnote{Department of People and Technology}\\
keld@ruc.dk
}

\maketitle

\begin{abstract}
Last-mile routing refers to the final step in a supply chain, delivering packages from a depot station to the homes of customers.
At the level of a single van driver, the task is a traveling salesman problem.
But the choice of route may be constrained by warehouse sorting operations, van-loading processes, driver preferences, and other considerations, rather than a straightforward minimization of tour length.
We propose a simple and efficient penalty-based local search algorithm for route optimization in the presence of such constraints, adopting a technique developed by Helsgaun to extend the Lin-Kernighan-Helsgaun algorithm for the traveling salesman problem to general vehicle routing problems.
We apply his technique to handle combinations of constraints obtained from an analysis of historical routing data, enforcing properties that are desired in high-quality solutions.
Our code is available under the open-source MIT license.
An earlier version of the code received the \$100,000 top prize in the Amazon Last Mile Routing Research Challenge organized in 2021.
\end{abstract}

\maketitle

\section{Introduction}

\footnotetext[1]{A technical report on this work was published online in \cite{lmrr2021}, as part of the submission to the Amazon Last Mile Routing Research Challenge. The current paper describes, in a more general context, an improved version of the local search algorithm adopted in the competition. }

One of the most visible examples of discrete optimization is the routing of delivery vans from doorstep to doorstep, bringing goods to individual customers in the {\em last mile} of a supply chain.
Vans are loaded, traveling salesman problem (TSP) tours computed, and drivers dispatched, all with the aim to minimize the time and cost of package delivery, while keeping drivers safe.

This general application is widely studied in the operations research literature, starting with a classic paper by \cite{dantzig1959}.
Particular focus has been given to variants of the capacitated vehicle routing problem (CVRP), where tour finding is combined with the task of assigning packages to vans.
Surveys of CVRP work can be found in \cite{laporte1992}, \cite{tv2002}, \cite{grw2008}, \cite{cruz2015}, and \cite{uchoa2017}.

In our work, we consider the last-mile problem at the level of a single driver, after warehouse operations have made a packages-to-vans assignment.
Although the typical target of 100 to 200 stops per van is well within the range of exact optimization by TSP solvers, the choice of routes may be subject to a variety of constraints not typically captured by the standard objective of minimizing tour length.
Some of these constraints, such as restrictions on the time-of-delivery for packages (\cite{christofides1981}), appear frequently and have been well studied in the literature.
Others can be specific to a particular distribution network, such as constraints imposed by sorting operations inside a warehouse or by the adopted van-loading process.
On top of these, additional constraints may arise due to driver preferences, regarding traffic, parking, and other factors.

To handle the wide variety of last-mile-routing problem types that can arise, we propose a local search algorithm, using penalties to enforce (or nearly enforce) all constraints, while aiming to minimize tour length.
The algorithm, based on the LKH-3 code by \cite{hel2017}, is very flexible, allowing the computer implementation to be quickly modified to handle new applications.

We demonstrate our search algorithm using data from the Amazon Last Mile Routing Research Challenge held in 2021 ({\url{https://routingchallenge.mit.edu}}).
In this context, constraints were obtained from an analysis of historical routes followed by van drivers.
Our entry into the competition received the \$100,000 top prize.

\subsection{Structure of this paper}
In Section~\ref{section_data} we describe the routing data provided by Amazon to competitors in the Last Mile Routing Research Challenge.
Our local search algorithm supporting disjunctive constraints is presented in Section~\ref{section_search} in a general form that is applicable for other similar TSP variants. This section includes  computational results for minimizing the tour length of the Amazon instances.
In Section~\ref{section_constraints} we describe the process used to extract routing constraints from Amazon's historical data so that our computed tours become similar to the actually driven tours.
We also show how these constraints are handled by the local search algorithm.  The constraints are presented in the order in which we observed and developed them during the competition.
Results from the competition are discussed in Section~\ref{section_contest},
together with a comment on the potential use of the
techniques in practice. Concluding remarks are given in Section~\ref{section_conclusions}.

\section{Amazon routing data}
\label{section_data}

The Last Mile Routing Research Challenge, run by Amazon and MIT, attracted 2,285 participants.
On top of the announced \$175,000 in prize money, a major draw was the opportunity to study routing data sampled from Amazon's last-mile network.
Indeed, the competition got underway with a massive release of 6,112 routes, described in \cite{amz2022}.
This is far and away the largest collection of real-world TSP instances ever made available to researchers.

Each of the 6,112 routes originates at one of 17 depots (called stations in the competition) in the United States and consists of a collection of dropoff stops to be visited by a driver on a specified day.
The number of stops ranges from a low of 32 to a high of 237.
The mean number of stops is 148.
For each stop to be visited, the data include the travel time in seconds to all other stops on the route, a zone ID, delivery time window, latitude and longitude, information on the sizes of packages to be delivered, and more.

The stop-to-stop travel times are asymmetric, that is, due to traffic flow, the estimated time to travel from stop $x$ to stop $y$ may not be the same as the time to travel from $y$ back to $x$.
The task of finding the quickest route for the driver is an instance of the asymmetric traveling salesman problem (ATSP).
Although direct optimization of the routes is not the focus of the challenge, the Amazon library of 6,112 practical instances of the ATSP will be a valuable resource for the operations research community.

Importantly, for the 6,112  routes, the data also include the actual sequence of stops in the tour that was followed by the van driver.
Each of the driver tours is classified as ``High" quality (2,718 instances), ``Medium" quality (3,292 instances), or ``Low" quality (102 instances) by Amazon for undisclosed reasons.
Rather than the standard TSP objective of minimizing total travel time, the competition asks for tours that closely match the ``High'' quality routes taken by experienced drivers.
We discuss this further in Section~\ref{section_constraints}.

\begin{figure}[htb]
  \centering
  \includegraphics[width=\columnwidth]{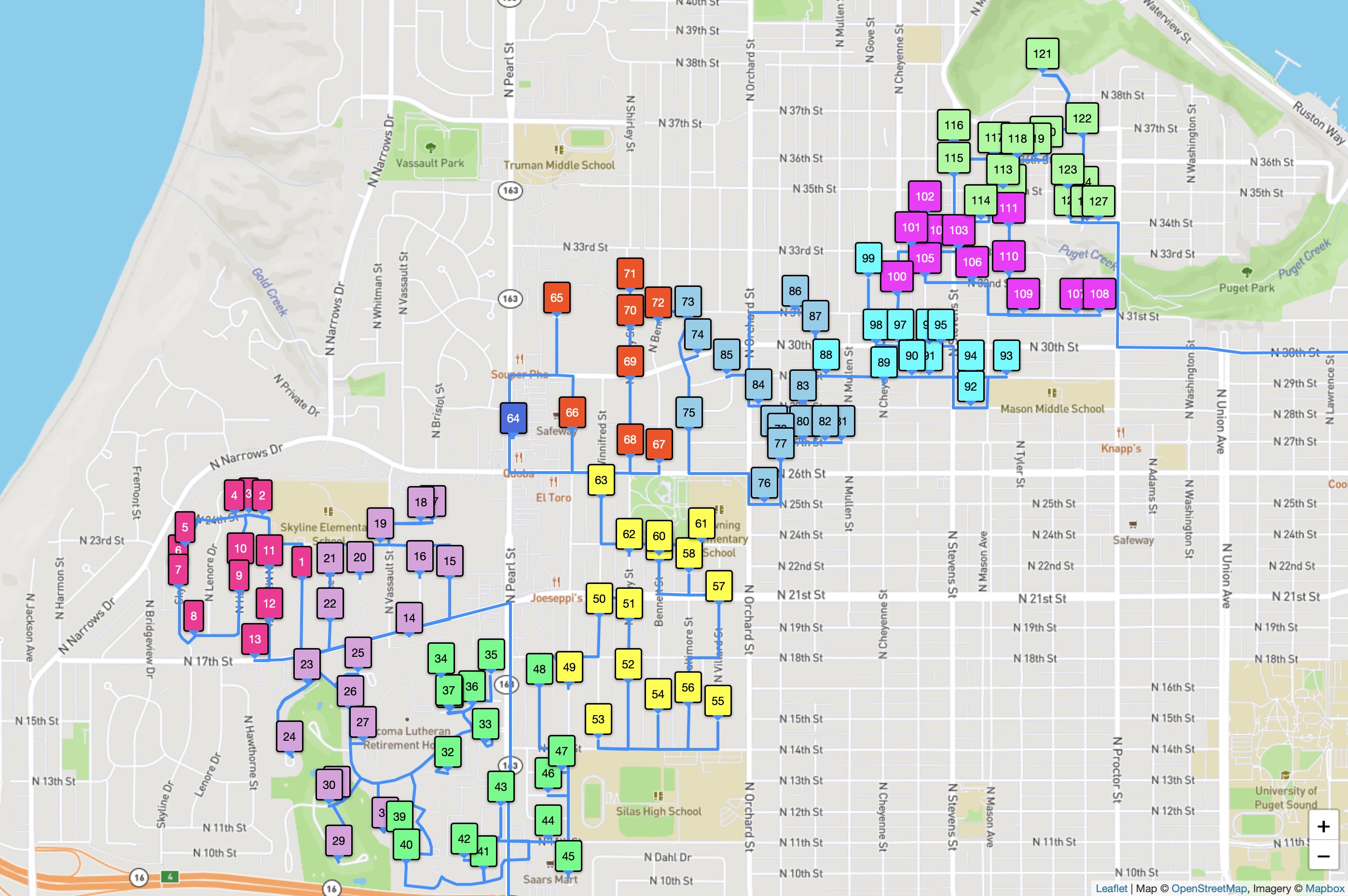}
  \caption{Amazon driver tour for route amz0002.}
  \label{figdriver}
\end{figure}

An illustration of one of the driver tours is given in Figure~\ref{figdriver}.
The labels indicate the order the stops were visited by the driver; the lines leaving the bottom and right-hand-side of the image indicate the paths from and to the depot, respectively.
Stops are colored by zone ID.
(The stop-to-stop paths were computed with the Open Source Routing Machine, {\url{http://project-osrm.org}}. Additional plots can be found on the page {\url{https://www.math.uwaterloo.ca/tsp/amz/maps.html}}.)

A portion of the driver tours include return visits to a stop due to an undelivered package during an initial visit.
The routes in the final scoring phase of the challenge are limited to those with no such undelivered packages (so no return visits to stops).
There are 1,107 such routes  that also receive the ``High'' quality rating.
This High+Delivered collection of instances will be the test set for our discussion.

As a benchmark, we computed optimal ATSP tours for these 1,107 test instances.
The Concorde code (\cite{abcc2006}) solves symmetric instances of the TSP, so we applied a standard transformation of the ATSP to the TSP, splitting each stop into a pair of nodes, one for incoming edges and the other for outgoing edges (see \cite{jt1983}).
The TSP instances thus obtained were solved with Concorde in an average of 142 seconds of computation time on a single core of a Linux workstation, equipped with an Intel Xeon Gold 6238 CPU @ 2.10GHz processor.
This is not an efficient way to solve instances of the ATSP, but it allowed us to easily obtain optimal tours for the test set. See Table~\ref{tab:atsp-driver-tours}.

\begin{table}[htb]
\caption{Total travel time of optimum ATSP vs. driver tours for the 1,107 selected test instances.}\label{tab:atsp-driver-tours}
\begin{displaymath}
    \begin{tabular}{cc} \hline
    Solution Set  & Mean Travel Time \\ \hline
    ATSP Tours    & 10853.3s         \\
    Driver Tours  & 12250.0s         \\ \hline
    \end{tabular}
\end{displaymath}
\end{table}

The driver tours are on average 12.9\% longer than optimal tours, highlighting the fact that the routes taken by drivers (based on tours proposed by Amazon) are likely taking into account constraints that are not captured by travel-time minimization alone.
Thus, in our approach in the competition, we use travel time only as a starting point,  mimicking driver deviation by adding constraints to the ATSP model.

\section{Constrained local search}
\label{section_search}

The general local search paradigm is simple: starting with any candidate solution to a problem,  look for small alterations that can possibly lead to a better solution.
The idea dates back to work on the TSP, starting with the 2-opt heuristic by \cite{flo1956}.
In Flood's algorithm, the small alteration consists of removing two edges in a TSP tour (considering a tour as a circuit in a graph) and reconnecting the resulting paths into a new tour, by adding two other edges.
This was quickly extended to 3-opt, removing and adding 3 edges at a time, by \cite{boc1958} and \cite{cro1958}, and later to a variable $k$-opt algorithm by \cite{lk1973}.

Building on the success of local search, extensions of the paradigm are discussed in surveys by \cite{lms2003} (iterated local search), \cite{hs2005} (stochastic local search), and \cite{avta2018} (guided local search), and in the book \cite{al2003}, covering a range of applications in discrete optimization.

Of particular interest in our work is the use of local search in the area of constraint programming, where penalties are used to guide the heuristic towards feasible solutions; see \cite{ht2006} for a survey.
Ideas from this constraint-programming work were adopted repeatedly to handle time window constrained routing problems, e.g.\  in \cite{lb2010},  \cite{nod2010},  or  \cite{vidalTimeWindows2013}.

\cite{hel2017} incorporated penalty-based search into his LKH-3 code for constrained traveling salesman and vehicle routing problems, greatly extending the reach of his well-known TSP code LKH (\cite{hel2000}).
We use in our work the simple and efficient method deployed in LKH-3, to measure simultaneously the length of the route and the penalties incurred when constraints are violated.

Our tour-finding heuristic is called LKH-AMZ, acknowledging its origins in LKH-3 and in the Amazon research challenge.
LKH-AMZ streamlines LKH-3 by implementing a restricted set of $k$-opt moves, suitable for the short computation times called for in the challenge.
At the same time, LKH-AMZ expands the family of constraints that can be handled, giving added flexibility.
We take advantage of this in our development of constraints to model driver preferences, described in Section~\ref{section_constraints}.

\subsection{Penalties}
\label{section_penalties}

To drive the heuristic, we assume the set of constraints is evaluated by a nonnegative-valued function, $pen(\cdot)$, taking as input a tour $T$ and returning 0 if and only if all constraints are satisfied.
The function can be as simple as the number of violated constraints, or a sum of component functions indicating for each constraint a measure of its violation.
We give examples in Section~\ref{section_constraints}.

Along with individual constraints, LKH-AMZ can support logical combinations, given in disjunctive normal form. An example is
to require
$C_1 \vee C_2 \vee C_3$, where $C_1, C_2, C_3$ are individual constraints.
In our application, we include disjunctions, evaluated as the minimum penalty among the component constraints.
(We note that a TSP application with disjunctive
constraints was recently introduced to tackle the cable tree wiring
problem in \cite{Koehler2021}.)

For a tour $T$, we let $len(T)$ denote its travel length.
Rather than attempting to combine $len(T)$ and $pen(T)$ into a single measure, we instead evaluate $T$ based on the pair of values.
We say a tour $\hat{T}$ {\em improves} $T$ if $len(\hat{T}) < len(T)$ and $pen(\hat{T}) \leq pen(T)$, and we seek improving $k$-opt moves.
The asymmetry in the use of $len(\cdot)$ and $pen(\cdot)$ is due to the mechanism we use in the search, considering only edge sets that strictly decrease the length of the current tour.

\subsection{Tour representation}

\begin{figure}[htb]
  \centering
  \includegraphics[width=3.0in]{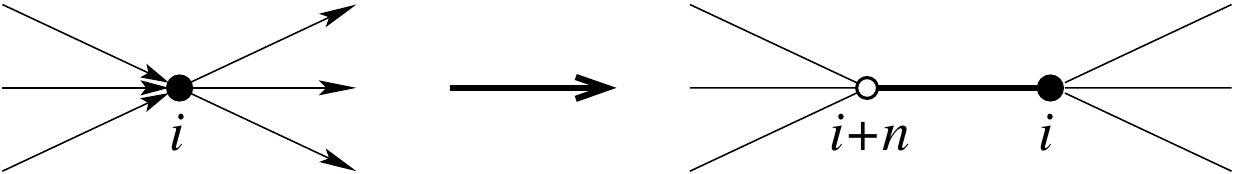}
  \caption{ATSP to TSP transformation.}
  \label{figtsp2atsp}
\end{figure}

To take advantage of data structures and algorithms specific to the symmetric TSP, LKH-AMZ adopts internally an ATSP to TSP transformation.
This increases the number of nodes $n$ in the ATSP instance to $2n$ in the TSP instance, by creating a dummy (``incoming'') node $i+n$ for each node $i$ in the ATSP instance.
Associated with each dummy node, we have a fixed edge $(i,i+n)$ that is required to be in the TSP tour; the fixed edge has 0 travel length.
For each original ATSP node $j$, other than $i$, the incoming edge $(j,i)$ is replaced by an edge $(j,i+n)$ having the same travel length as $(j,i)$.
See Figure~\ref{figtsp2atsp}.

In the transformed problem, a TSP tour can be oriented in one of two ways.
To simplify the code, LKH-AMZ organizes its tours so that every original node follows its dummy node in the selected orientation.

The heart of the LKH-AMZ algorithm is a process to search for an improving $k$-opt move.
To support the search, the current tour $T$ is stored both as a doubly-linked list, giving immediate access to the predecessor and successor of each node in the tour orientation, and as a $rank$ array, listing the position of each node in the tour order.
The array allows the code to quickly test if a node $b$ lies between two other nodes $a$ and $c$ in the tour order, that is, $rank(a) < rank(b) < rank(c)$.
These fast queries give an $O(1)$ computation to list the $2k$ edges involved in a proposed $k$-opt move, described in Section~\ref{subsection_kopt}.

In the search algorithm, a move is considered only if the length of the $k$ added edges is strictly less than the length of the $k$ deleted edges.
In this case, a doubly-linked-list representation of the new tour $T'$, produced by the $k$-opt move, is obtained by adjusting $2k$ pointers in the representation of $T$.
Our penalty function is designed to work directly with this linked-list representation, not needing the $rank$ array.
Thus, only if $pen(T') \leq pen(T)$ do we actually complete the $k$-opt move by updating the array, an $O(n)$ computation.

In our study, the penalty test rejects over 90\% of proposed moves, justifying the use of the pair of data structures, despite the computationally expensive update.
(We thank Francesco Cavaliere for making this observation, leading us to replace the two-level tree data structure (see \cite{fjmo1995}) we adopted originally in LKH-AMZ.)

\subsection{Candidate edges}

LKH-AMZ restricts the search for $k$-opt moves by considering for possible inclusion in the tour only a subset of the edges, referred to as {\em candidate edges}.
The algorithm takes as input a parameter {\tt MAX\_CANDIDATES} that bounds the number of incident candidate edges at each node; by default the parameter is set to 6.

To select candidate  edges, we begin by employing the iterative method of \cite{heldkarp1971} to compute a 1-tree relaxation $S$ of the TSP (as adopted in \cite{hel2000}).
For each edge $(i,j)$, we set $\alpha(i,j)$ equal to its reduced cost with respect to $S$, that is, the amount the length of $S$ would increase if $(i,j)$ is required to be in the 1-tree.
Then, for each node $i$, we select the {\tt MAX\_CANDIDATES} edges meeting $i$ having the smallest $\alpha(\cdot)$ values.

\subsection{Improving \textit{k}-opt moves}
\label{subsection_kopt}

In our search algorithm, we consider single $k$-opt moves that improve the current tour, rather than employing the variable $k$-opt method of \cite{lk1973}, where small (possibly non-improving) moves are chained together.
In this sense, the ``LK'' in the algorithm name is misleading.
It is not a Lin-Kernighan algorithm, but rather a version of $k$-opt.

Furthermore, to reduce the number of calls to the potentially time-consuming penalty function, we consider only moves with $k=3$ and $k=4$.
Note that 2-opt moves are not permitted, since they require a reversal of a subsequence of tour edges, violating the fixed orientation that keeps each original node as a successor of its corresponding dummy node.

Given a tour $T$, and a non-fixed edge $(t_1, t_2)$ in $T$, we attempt to find an improving move that deletes $(t_1, t_2)$.
Again, since every original node must follow its dummy node, a move is not allowed to reverse any segment of the tour.
With this restriction, we can consider only two possible types of moves.

\begin{figure}[htb]
  \centering
  \subfloat[Choice of points.]{\includegraphics[width=.25\columnwidth]{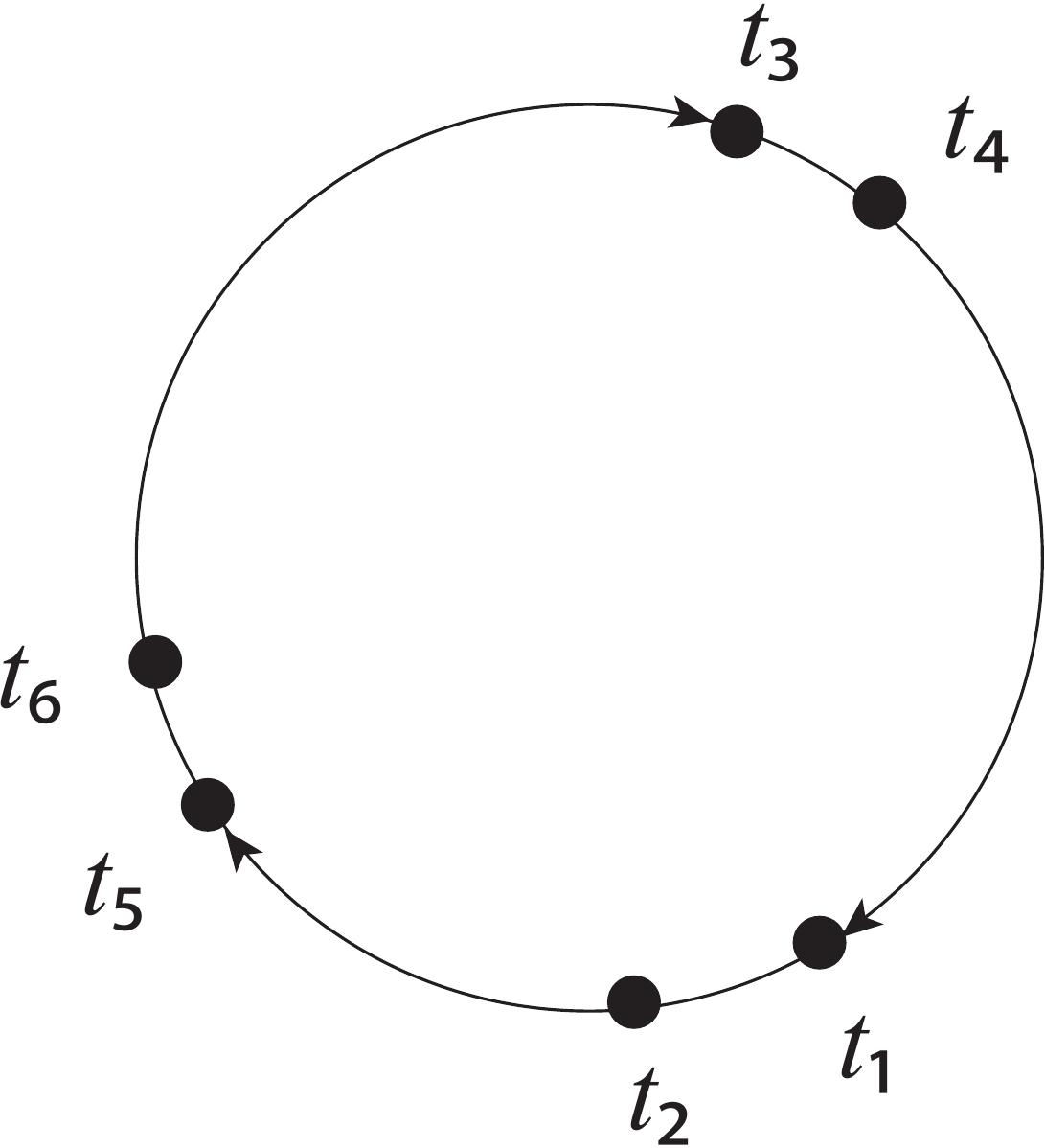}\label{fig3pts}}
  \hskip 4em
  \subfloat[Result of 3-opt move.]{\includegraphics[width=.25\columnwidth]{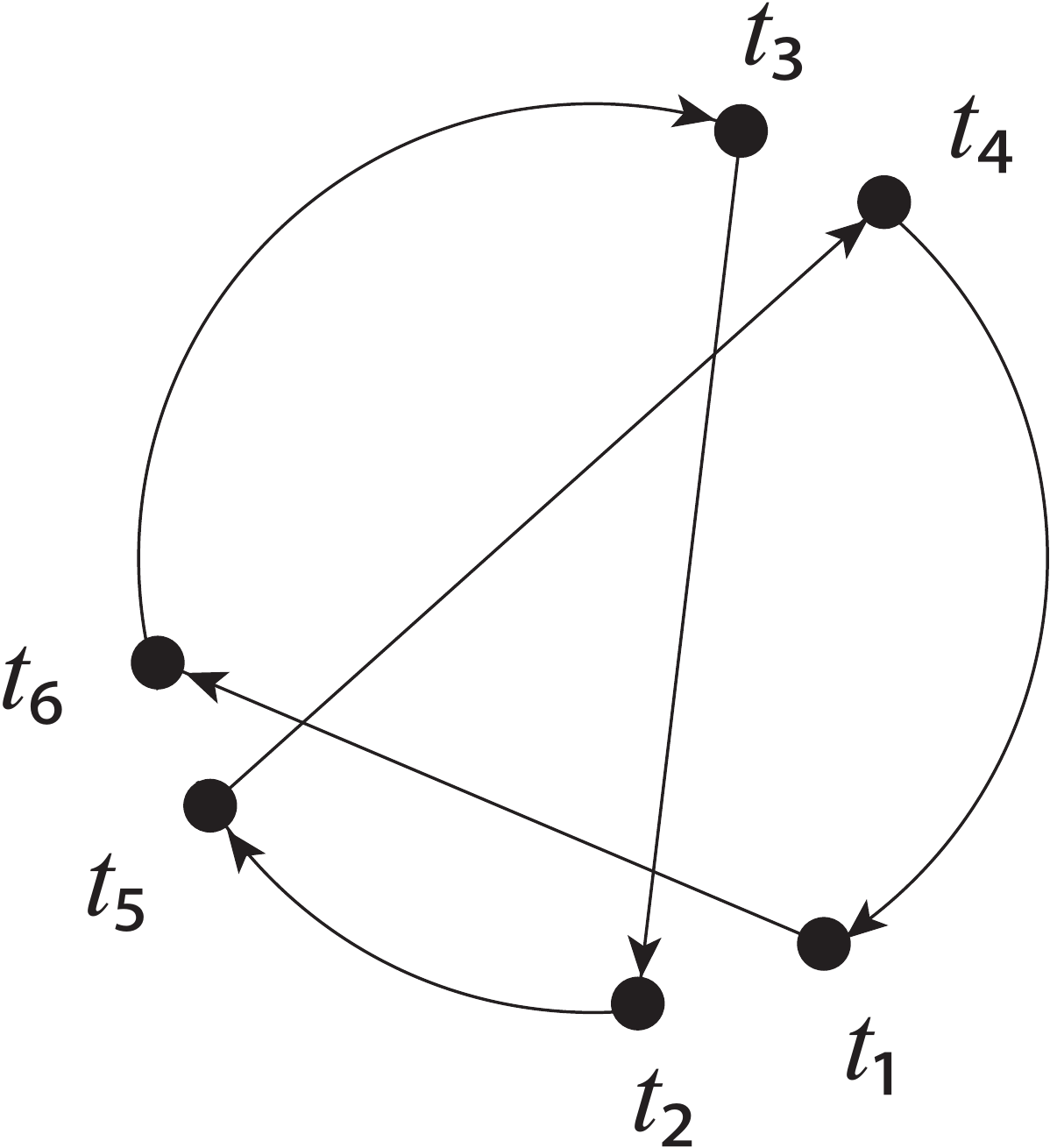}\label{fig3opt}}
  \caption{Special 3-opt move.}
  \label{fig3}
\end{figure}

The first of these moves is illustrated in Figure \ref{fig3}, involving nodes $t_1, t_2, t_3, t_4, t_5$, and $t_6$.
Note that node $t_5$ must lie between $t_2$ and $t_3$ in the tour's orientation.
Figure \ref{fig3opt} shows the result of the 3-opt move.
Observe that no segments have been reversed.

To find such a 3-opt move, we consider each candidate edge $(t_2,t_3)$ for node $t_2$, such that the travel length of $(t_2,t_3)$ is strictly less than the length of edge $(t_1,t_2)$.
This restriction on the choice of $(t_2,t_3)$ is known as the {\em positive gain criterion}, where candidate edges are chosen only if the cumulative gain (in reduced travel length) from the proposed set of edge exchanges is positive.
Letting  $t_4$ denote the successor of $t_3$ in the tour, we consider each of $t_4$'s  candidate edges $(t_4,t_5)$ that satisfy the positive gain criterion, that is,
$$length(t_2,t_3) + length(t_4,t_5) < length(t_1,t_2) + length(t_3,t_4).$$
We then let $t_6$ denote the successor of $t_5$.
We process the move only if the 6 nodes follow the pattern in Figure \ref{fig3pts}.

\begin{figure}[htb]
  \centering
  \subfloat[Choice of points.]{\includegraphics[width=.25\columnwidth]{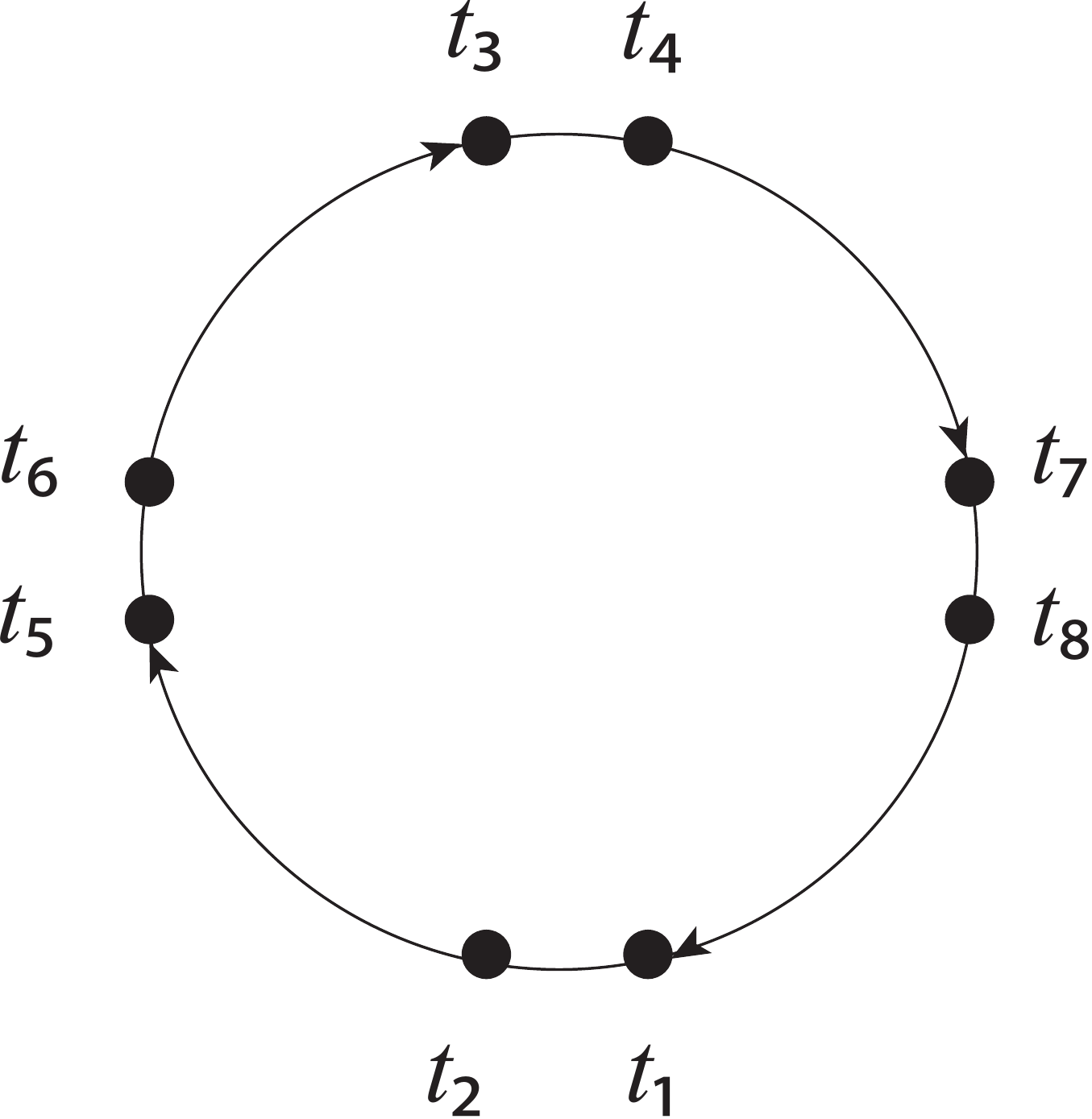}\label{fig4pts}}
  \hskip 4em
  \subfloat[Result of 4-opt move.]{\includegraphics[width=.25\columnwidth]{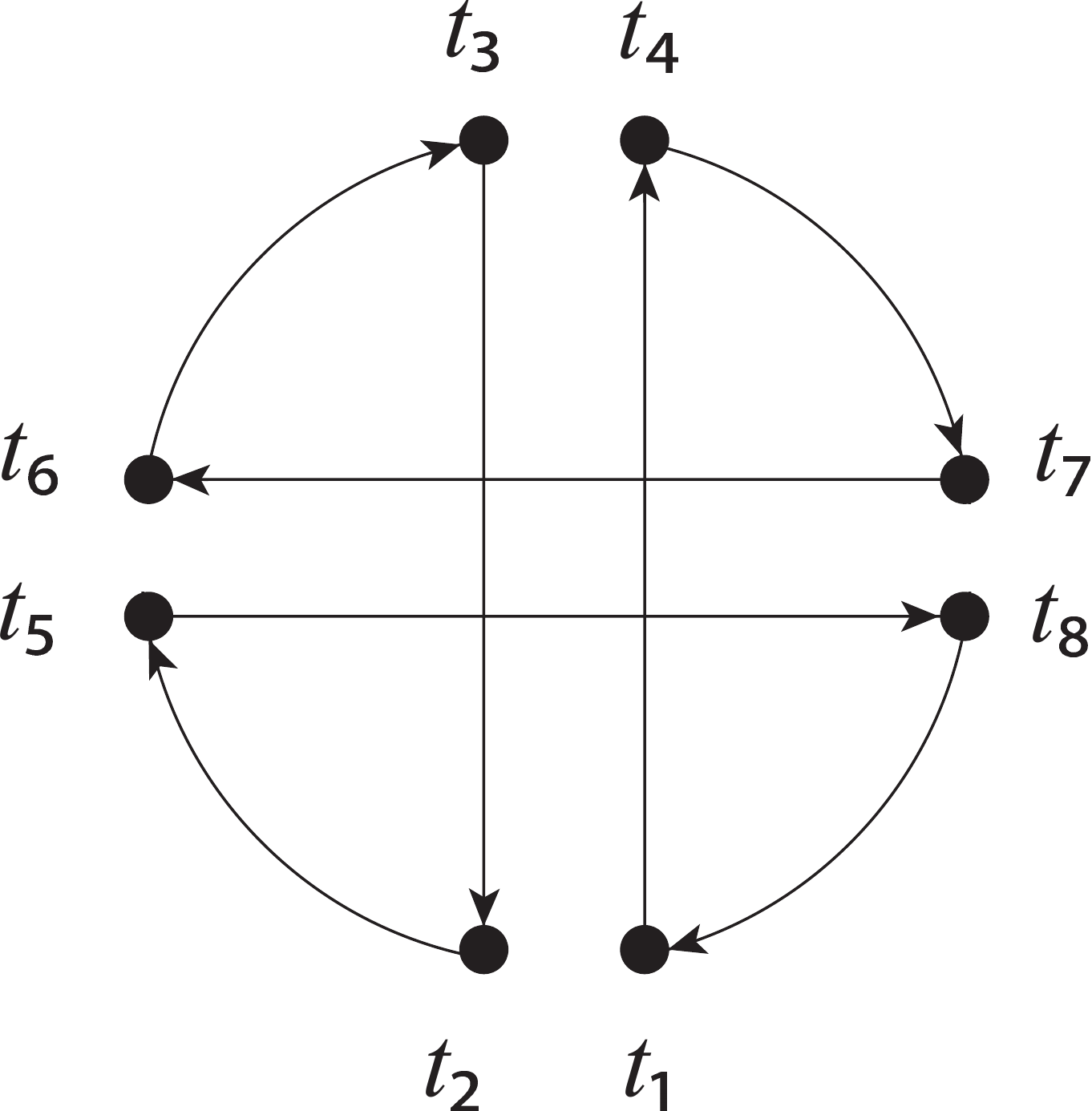}\label{fig4opt}}
  \caption{Special 4-opt move.}
  \label{fig4}
\end{figure}

The second of the moves is illustrated in Figure \ref{fig4}, with the choice of nodes $t_1, t_2, t_3, t_4, t_5, t_6, t_7$, and $t_8$.
Note that again $t_5$ must lie between $t_2$ and $t_3$, and $t_7$ must be between $t_4$ and $t_1$.
To locate such a 4-opt move, we use $t_2$'s candidate edges and the positive gain criterion to select $t_3$ and its successor $t_4$.
Now, to keep the time complexity low, rather than proceeding with $t_4$'s candidate edges, only four possible combinations of $(t_5, t_6)$ and $(t_7, t_8)$ are tried: those nearest to $t_1$, $t_2$ and $t_3$, $t_4$, respectively.

\subsection{Iterated local search}
\label{sec_iterated_ls}

A {\em trial} in LKH-AMZ is a sequence of 3-opt and 4-opt moves, where the search for these moves is terminated according to rules that trade off the speed of the computation with the likelihood of finding further improvements.
(See \cite{abcc2006} for a detailed discussion in the context of the pure TSP.)
In each {\em run} of the code, $8n$ trials are executed.

The trials are organized as an iterated local search.
The idea is to allow the code to quickly sample local minima in the neighborhood of the best tour $T^*$ that has been found thus far.
The initial trial begins with a pseudo-random tour.
In each subsequent trial, the current $T^*$  is ``kicked'' (perturbed) by performing a 4-opt move of the type indicated in Figure~\ref{fig4}; the 4-opt kick is selected using the Rohe long-edge rule with a 50-step random walk, described in \cite{acr2003}. In the Rohe rule, the first step of the random walk aims to remove a long edge. To this end it chooses a stop $v$ from a small randomly selected set that maximizes $length(v,next(v)) - length(v,near(v))$, where $next(v)$ is the next stop after $v$ in the tour and $near(v)$ is its nearest neighbor.
The kicked tour is re-optimized  via a sequence of 3-opt and 4-opt moves, producing a tour $T$.
If the travel length of $T$ is not less than the travel length of $T^*$, then we attempt to replace segments of $T$ by segments of $T^*$ in a process known as {\em iterative partial transcription} (IPT), proposed by \cite{moebius1999}.
If tour $T$ (possibly updated via the IPT process) improves $T^*$, then we replace $T^*$ by $T$.

\subsection{Multiple runs}
\label{sec_multiple_runs}
Until a specified time limit is reached, we carry out multiple runs of the full iterated local search.
After each run, we evaluate the resulting tour $T$ by the function $$v(T) = PenaltyMultiplier * pen(T) + len(T)$$ where the scaling factor $PenaltyMultiplier$ (set to 1500 by default in the Amazon tests) allows us to balance the penalty and length factors.
We record $T$ as a new best tour if $v(T)$ is less than the $v(\cdot)$ value of our current best tour.

\subsection{Computational results on ATSP instances}
\label{subsection_atspcomp}

The LKH-AMZ heuristic is optimized to handle the possibly large number of constraints that can arise in practical applications.
Nonetheless, as shown in Table~\ref{table_atsp}, the search method achieves near-optimal tours on the 1,107 High+Delivered unconstrained instances in the Amazon data set.
\begin{table}[htb]
\caption{LKH-AMZ runs on 1,107 ATSP instances}
\label{table_atsp}
\begin{center}
\begin{tabular}{c|cc|cc} \hline
    \multicolumn{1}{c}{ } & \multicolumn{2}{c}{Intel Xeon Gold} & \multicolumn{2}{c}{Apple Silicon M1}\\ \hline
Run Time & Opt Ratio & Opt Tours & Opt Ratio & Opt Tours \\ \hline
    1s   & 1.000101  & 998       & 1.000039  & 1058      \\
    2s   & 1.000023  & 1070      & 1.000007  & 1095      \\
    5s   & 1.000005  & 1100      & 1.000002  & 1103      \\
    10s  & 1.000002  & 1104      & 1.000001  & 1105      \\
    20s  & 1.000000  & 1106      & 1.000000  & 1107      \\
    30s  & 1.000000  & 1107      &           &           \\ \hline
\end{tabular}
\end{center}
\end{table}

We report results from runs on a Linux server equipped with an Intel Xeon Gold 5218 processor (2.30GHz) and on a 2020 MacBook Air equipped with the Apple Silicon M1 processor.
The M1 processor gives an impressive performance, producing results on the MacBook that are significantly better than the server-class Intel processor.
Using 1 second of running time per instance, the code on the MacBook found tours of average length within a factor of 1.000039 of optimal (that is, on average less than 0.5 seconds longer than an optimal tour).
With 20 seconds of running time, it found optimal tours in all 1,107 instances.

\section{Adding constraints to predict driver tours}
\label{section_constraints}

In the evaluation phase of the Last Mile Routing Research Challenge, submitted codes were scored on a set of 3,072 routes that were kept hidden from participants.
The new routes were expected to have properties similar to the set of 1,107 High+Delivered instances described in Section~\ref{section_data}.
Thus, our strategy in the competition was to use the training data to discover constraints satisfied by high quality tours, allowing us to employ LKH-AMZ as an optimization engine.

As we stated earlier, rather than minimizing tour distance, teams in the competition were asked to match, as best they can, the sequences of stops in actual routes followed by van drivers.
From a practical viewpoint, this seems at first glance a surprising choice.
If drivers can find good routes, why are new optimization methods needed?
A short answer is that drivers are not building their routes from scratch.
Indeed, the announcement of the challenge on the Amazon Science blog (\cite{amazon2021}) indicates drivers are presented with computed tours that they adjust on-the-fly.
\begin{quote}
``Drivers, however, frequently deviate from those computed routes. Drivers carry information about which roads are hard to navigate, when traffic is bad, when and where they can easily find parking, which stops can be conveniently served together, and many other factors that existing optimization models don't capture.''
\end{quote}
These computed tours were not made available, but the scoring metric directed challenge competitors towards global features (likely due to Amazon algorithms), together with local adjustments (made by drivers, based on their preferences).
The hoped-for outcome is an algorithm to quickly compute efficient tours that drivers are happy to follow.

In the remainder of this section, we will briefly describe the scoring method of the contest. Then, we derive additional constraints for the ATSP in the order in which we observed structures in the benchmark data during the competition.
For each new constraint type, we show how it reduces the score of the planned tours.
  Most obvious are the zone IDs that lead to a clustered ATSP (Section~\ref{subsection_zone_ids}), while time windows did not play a central role for us (Section~\ref{subsection_time_windows}). An important ingredient to improve the score was to predict the order of zones in a tour.
In Section~\ref{subsection_order} we show how to model them as constraints in our local search.
Then, in Section~\ref{subsection_pred} we show how to learn zone orders from the training data and in Section~\ref{subsection_cluster}
how to guess hierarchical orders from the structure of the zone IDs.
We conclude with a brief discussion on the complexity of evaluating the penalties for all constraints in Section~\ref{section_complexity_penalties}.

\subsection{Scoring}

The scoring metric used in the competition is complex, combining two ways to measure the difference between two tours, one involving only the positions of the stops in the tours and the other also taking into account the travel time between stops.
The two measures are multiplied and scaled, to produce a single numerical score, the lower the better, with a perfect match of the driver route receiving 0.0.
Details of the score computation can be found at the challenge GitHub site {\url{https://github.com/MIT-CAVE/rc-cli/tree/main/scoring}}.
The overall score for a competing team is the average of their tour scores, taken over the full hidden data set.

To indicate the range of outcomes, random tours on the High+Delivered instances received a score of 0.91545, averaged over 10 trials.
Much better are the optimal ATSP tours, attaining a score of 0.07030.
Our target is to find routing constraints that will drive this score towards 0.0.

\subsection{Zone IDs}
\label{subsection_zone_ids}

\begin{figure}[htb]
  \centering
  \includegraphics[width=\columnwidth]{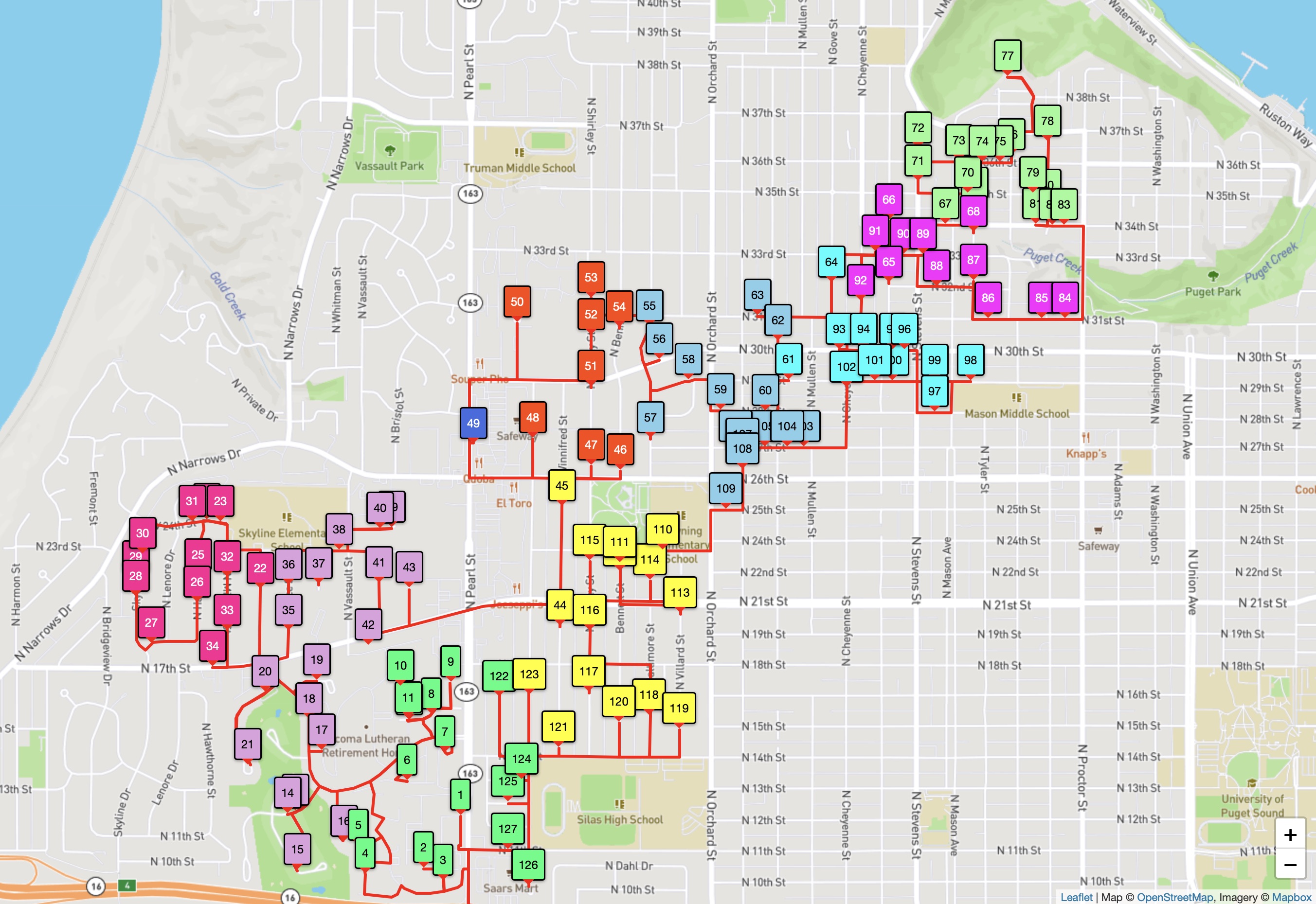}
  \caption{Optimal ATSP tour for route amz0002.}
  \label{figatsp}
\end{figure}

To begin, compare the driver tour displayed in Figure~\ref{figdriver} with the optimal ATSP tour for the same instance displayed in Figure~\ref{figatsp}.
The ATSP tour arrives from the depot and returns to the depot via the path displayed at the bottom of the map.
Like in the earlier figure, the order of stops in the tour is indicated by the numbers on the labels.
Also, in both figures, the coloring (shading) of the labels indicate stops that share a common zone ID.

The tours are similar (hence the 0.07030 score), but notice that the driver tour visits consecutively all stops having the same zone ID, that is, each label color corresponds to a cluster in the tour order.
The ATSP tour does not share this property.
This is the first feature we will exploit to improve the score of our tours.

In the High+Delivered data set, there are an average of 20 distinct zone IDs represented in each route, giving an average of approximately 7 stops per zone.
The zones partition the set of nodes in the ATSP instance, where we create a special zone containing only the depot.
(In cases where a stop in the training data is missing a zone ID, we assign to it the ID of the closest stop (in Euclidean distance) associated with the same depot.)
We use this zone partition to create an instance of the {\em clustered ATSP}, where the stops in each zone are required to appear consecutively in any candidate tour.

\begin{figure}[htb]
    \centering
  \includegraphics[width=\columnwidth]{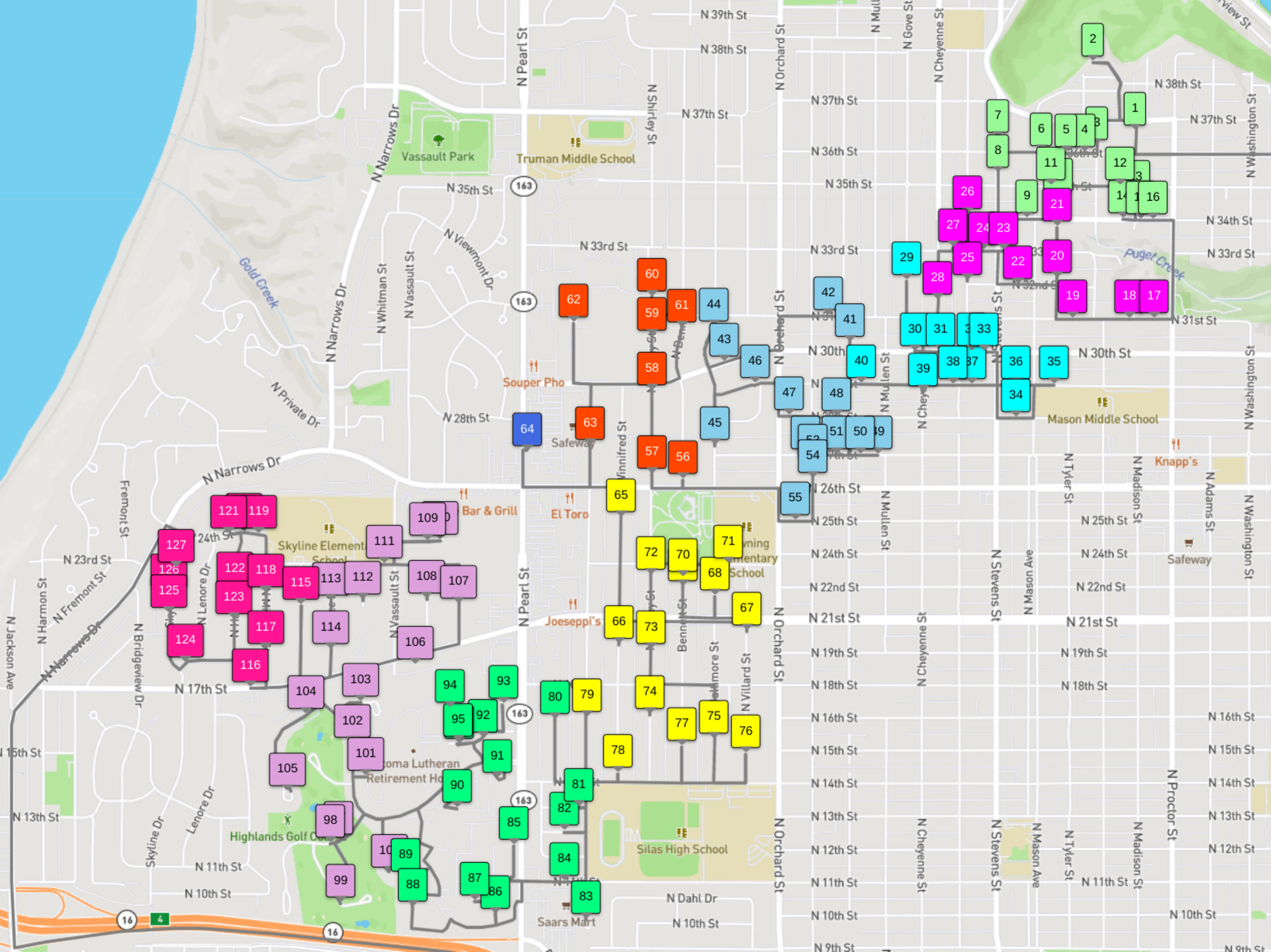}
  \caption{Clustered ATSP tour for route amz0002.}
  \label{figcluster}
\end{figure}

An instance of the clustered TSP can be transformed to an ATSP instance by adding a large constant $M$ to the travel time of any edge joining stops in distinct clusters.
Applying this modification to the High+Delivered test set, we again used the ATSP to TSP transformation and solved the resulting instances with Concorde.
\begin{table}[htb]
\caption{Comparison of ATSP, clustered, and driver tours.}
\label{tab:atsp-cluster-driver}
\begin{displaymath}
    \begin{tabular}{ccc} \hline
    Solution Set     & Mean Travel Time & Score   \\ \hline
    ATSP Tours       & 10853.3s         & 0.07030 \\
    Clustered Tours  & 11235.1s         & 0.04866 \\
    Driver Tours     & 12250.0s         & 0.00000 \\ \hline
    \end{tabular}
\end{displaymath}
\end{table}
The forced clustering increases the length of the tours by 3.5\% on average, but improves the score to 0.04866 as shown in Table~\ref{tab:atsp-cluster-driver}.
An image of a clustered tour can be seen in Figure~\ref{figcluster}.

The clustered ATSP instances give a second benchmark for the LKH-AMZ heuristic.
The constraints can be handled either by the transformation to the ATSP or by a penalty function based on the excess number of inter-cluster edges in the tour.
To be more precise, if there are $z$ clusters and for a tour $T$ we let $crossing(T)$ denote the number of times $T$ enters a new cluster, then
$pen(T) = \rho * (crossing(T) - z)$,
where $\rho$ is a constant.
The value of $pen(T)$ can be computed by a single pass through the tour order.

In Table~\ref{table_ctsp} we report results on both ways of implementing the constraints.
\begin{table}[htb]
\caption{LKH-AMZ runs on 1,107 clustered ATSP instances (M1 processor)}
\label{table_ctsp}
\begin{center}
\begin{tabular}{c|ccc|ccc} \hline
\multicolumn{1}{c}{ } & \multicolumn{3}{c}{ATSP Transformation} & \multicolumn{3}{c}{Cluster Penalties}\\ \hline
Run Time & Opt Ratio  & Opt Tours & Score & Opt Ratio & Opt Tours & Score \\ \hline
1s  & 1.000052 & 1051 & 0.04847 & 1.000302 & 931  & 0.04862 \\
2s  & 1.000010 & 1092 & 0.04846 & 1.000112 & 1021 & 0.04862 \\
5s  & 1.000002 & 1102 & 0.04850 & 1.000019 & 1085 & 0.04864 \\
10s & 1.000000 & 1106 & 0.04853 & 1.000005 & 1098 & 0.04856 \\
20s & 1.000000 & 1107 & 0.04853 & 1.000003 & 1104 & 0.04860 \\ \hline
\end{tabular}
\end{center}
\end{table}
The competition scores of the ATSP transformation and the cluster penalties are similar, and, in fact, slightly better than the score for the optimal clustered tours found by Concorde as shown in Table~\ref{tab:atsp-cluster-driver}.
But using the cluster transformation allowed LKH-AMZ to find shorter tours on average.
This can be explained, in part, by the better quality of the candidate edges produced when travel times incorporate important constraints such as cluster constraints.
We will therefore adopt the transformation in the computations discussed in the remainder of the paper.

\subsection{Time windows}
\label{subsection_time_windows}

The Amazon data include estimates of the service times at each stop, and, for a subset of the stops, time window targets for package delivery.
These targets can be handled as constraints, where a penalty is incurred if a tour $T$, including the service times, arrives at a stop outside  its time window.
In our implementation, $pen(T)$ is computed by a single pass through the tour, summing the number of late  seconds.
Early arrivals at a stop lead to a waiting time until the time window starts. This way, early arrivals do not lead to a violation. We do not try to reduce the number of late violations by allowing early violations.

The computational results on the High+Delivered data set, reported in Table~\ref{table_window}, reflect the fact that only a small number of stops have non-trivial time windows.
The ``ATSP Transformation without TW opt.'' runs are essentially equivalent to those in Table~\ref{table_ctsp}, ignoring time windows during the optimization, but evaluating the number of violations at the end.
  The ``ATSP Transformation with TW opt.'' runs   penalize late arrivals during the optimization.

\begin{table}[!thb]
\caption{LKH-AMZ scores on 1,107 clustered ATSP instances with time windows  (M1 processor)}
\label{table_window}
\begin{center}
  \begin{tabular}{c| c c  | c c } \hline
  \multicolumn{1}{c}{} & \multicolumn{2}{c}{ATSP Transformation without TW opt.} & \multicolumn{2}{c}{ATSP Transformation with TW opt.}\\ \hline
Run Time  &  Score & \#Violations & Score & \#Violations  \\ \hline
1s        &  0.04847 & 695 & 0.04843  & 0 \\
2s        &  0.04846 & 692 & 0.04848  & 0 \\
5s        &  0.04850 & 674 & 0.04856  & 0 \\
10s       &  0.04853 & 667 & 0.04861  &  0 \\
20s       &  0.04853 & 682 & 0.04861  &  0 \\ \hline
\end{tabular}
\end{center}
\end{table}

Ignoring time windows leaves at most 695 violations among the 163,177 stops in the given instances. Using penalization, no time window violations remain.
  The score  difference is well within the expected variability (compare the ``ATSP Transformation'' results in Table~\ref{table_ctsp}).
Notably, the Amazon driver tours also leave 647 time window violations on these instances. This is only slightly less than what our solutions ignoring them altogether achieve.
Given the negligible impact on the competition score, we do not include time window constraints in our further constraint generation.
At the end of Section~\ref{subsection_cluster}, we show that including time window constraints into our best known constraint set significantly degrades the results.
\subsection{Driver zone order}
\label{subsection_order}

The design of LKH-AMZ allows the code to handle constraints involving individual stops, such as the time window targets.
But the 6,112 training routes in the Amazon data include few examples where stops or streets are repeated.
This makes it difficult to extract meaningful constraints at this level.

However, the shortest clustered ATSP tour in Figure~\ref{figcluster} still differs significantly from the driver tour in Figure~\ref{figdriver}. While the driver serves the zones from the south-west towards north-east, the  shortest clustered ATSP tour
takes the opposite direction (the stop order is given by the stop labels).
We  focus our attention on constraints guessing the order by which the  zones are traversed in each tour.
As a first experiment, we present results for tours that are constrained to follow the same order of zones as in the driver tour.
(In cases where the driver tour splits a zone into two or more subsequences, we select the first occurrence when constructing the zone order.)
Of course, in a new planning phase we do not know the preferred zone order of the driver.
Learning and generating meaningful zone order constraints will be discussed later in  Sections~\ref{subsection_pred} and \ref{subsection_cluster}.

To enforce the zone ordering, we introduce two types of inter-zone constraints.
The first of these requires a pair of zones $(a,b)$ to be neighbors in the tour order, that is, if we let $visit(x)$ denote the position of zone $x$ in the tour (with the depot zone having position 0), then a {\em neighbor constraint} for $(a,b)$ requires $|visit(a) - visit(b)| = 1$.
Note that such a constraint can be handled by a travel-time transformation, adding a large constant $M$ to the travel time of all edges having one end in the union of zones $a$ and $b$ and the other end not in the union.

A second constraint for the zone pair $(a,b)$ is the stronger restriction that zone $a$ immediately precedes zone $b$ in the tour, that is, visit($a$) $=$ visit$(b)-1$.
This {\em path constraint} can also be handled by a travel-time transformation, first enforcing the neighbor constraint for $(a,b)$ and then adding the constant $M$ to the travel time of every edge directed between nodes in zone $b$ and nodes in zone $a$.

Alternatively, the two classes of constraints can be handled in the penalty function $pen(\cdot)$.
To do this efficiently (and to handle further zone-level constraints), we make a single pass through the tour order, computing the $visit(\cdot)$ value for each zone.
If a zone $z$ is split in the tour $T$, then we set $visit(z)$ by the position of its final subsequence.
With these values, we have a fast test for checking individual constraints, and we add to $pen(T)$ a constant $\rho_{neigh}$ for each violated neighbor constraint and a constant $\rho_{path}$ for each violated path constraint.

Let $z_1, \ldots, z_k$ be the zones in the order they appear in the driver tour. Adding a single path constraint $(z_1,z_2)$ to set the orientation of the tour and then  neighbor constraints $(z_i,z_{i+1})$ for $i = 2, \ldots , k-1$ can enforce this tour order.

\begin{table}[htb]
\caption{LKH-AMZ scores on driver zone-order instances (M1 processor)}
\label{table_fixed}
\begin{center}
\begin{tabular}{c|ccc|ccc} \hline
Run Time &  Transformed & Penalty \\ \hline
1s  & 0.00510 & 0.00738 \\
2s  & 0.00505 & 0.00748 \\
5s  & 0.00505 & 0.00778 \\ \hline
\end{tabular}
\end{center}
\end{table}
Results in Table~\ref{table_fixed} are presented for both ways of handling the constraints.
The scores show a dramatic improvement over those obtained when using the clustered ATSP instances (see Table~\ref{tab:atsp-cluster-driver}).
Also, note again an advantage for travel-time transformations over the penalties in a case where the entire constrained problem can be formulated as a transformed instance of the ATSP.

\subsection{Zone precedence constraints from training data}
\label{subsection_pred}

The competition scores reported in Table~\ref{table_fixed} are excellent, but we do not know beforehand the driver zone order when processing a new route $R$.
We therefore attempt to make use of the orders of the training routes, extracting constraints to guide LKH-AMZ towards a tour having similar characteristics.
As a preliminary step, we build a directed graph for every route, having a node for each visited zone (other than the depot) and a directed edge $(a,b)$ if the driver tour visits zone $a$ immediately before zone $b$.
Note that the graph may contain cycles, e.g. a driver might visit zones $a,b,c,d,e,f$ in the following sequence
  $$a\to b\to c\to a\to d\to e\to f\to e.$$
By contracting all strongly-connected components in the graph, we obtain a {\em component path} for the route (see \cite{tar1972}).
The nodes of the component path correspond to sets of zones.
In the above example we obtain the component path $\{a, b, c\} \to \{d\} \to \{f, e\}.$

We use component paths to construct {\em precedence constraints} for pairs of zones $(a,b)$ in route $R$, requiring $visit(a) < visit(b)$.
There are many ways to extract such constraints from the training instances.
In our code, when processing $R$, we select a single reference route $Q$ from the training set, such that $Q$ starts at the same depot as $R$ and maximizes the number of zones $Q$ has in common with $R$, multiplied by 2 for ``High'' routes, 1.5 for ``Medium'' routes, and 1 for ``Low'' routes.
From the component path of $Q$, we build a {\em pruned component path} by deleting any component that contains no zones in $R$ and adding edges to join the resulting path segments into a path (maintaining the same order as in the component path).
If the above example is our reference route $Q$ and our target route $R$ has zones $a,b,e, $ and $ g$, then the pruned component path is $\{a, b, c\} \to \{f, e\}$.

We create a precedence constraint (for route $R$) for each pair of zones $(a,b)$ such that $a$ and $b$ are in both $R$ and $Q$ and the strongly connected components $A\ni a$ and $B\ni b$ are joined by an edge $(A,B)$ in the pruned component path we obtained from $Q$.
In our example, we extract the precedence constraints
$visit(a) < visit(e)$ and $visit(b) < visit(e)$.
For zones $a$ and $b$ no precedence order can be derived, since they are visited cyclically in $Q$, while zone $g$ does not occur in $Q$.

As an alternative, we create additional precedence constraints from the edges in the transitive closure of the reference route's pruned component path.

Although there is not an efficient transformation to incorporate precedence constraints directly into the travel times, they can be handled with penalties in the same manner as the neighbor and path constraints discussed in Section~\ref{subsection_order}.

In our computational tests, reported in Table~\ref{table_pred}, we use the full 6,112-instance Amazon set for training.
\begin{table}[htb]
\caption{Scores of LKH-AMZ runs with precedence constraints (M1 processor)}
\label{table_pred}
\begin{center}
\begin{tabular}{c|cc} \hline
Run Time & Component Path & Transitive Closure \\ \hline
1s       & 0.03322        & 0.03171 \\
2s       & 0.03303        & 0.03169 \\
5s       & 0.03316        & 0.03169 \\ \hline
\end{tabular}
\end{center}
\end{table}
When processing a route $R$, we remove it from the training set, to avoid having the reference route $Q$ chosen as $R$ itself.
The table reports scores for adding to the basic clustered ATSP model either the precedence constraints from only the component path, or the much larger set of constraints from the transitive closure.
The results are a significant improvement over the 0.048xx scores obtained with the basic model.

We remark that several attempts to aggregate precedence constraints from multiple reference routes did not improve the single-reference results reported in Table~\ref{table_pred}.

\subsection{Cluster rules from zone IDs}
\label{subsection_cluster}
The Amazon zone IDs assigned to stops have the form $\Gamma$-$x.y\Delta$, where $\Gamma$ and $\Delta$ are capital letters and $x$ and $y$ are integers.
In Table~\ref{zone_table} we list these IDs in the order they appear in the driver tours for each of the first four High+Delivered routes.
\begin{table}[htb]
\caption{Sequences of zone IDs in driver tours}
{
\label{zone_table}
\begin{center}
\begin{tabular}{|l|l|l|l|}\hline
amz0002& amz0003& amz0012& amz0019 \\
\hline
Depot&  Depot&    Depot&  Depot \\
A-2.2E& B-11.1G& A-2.1E& M-10.3D \\
A-2.1E& B-11.1H& A-2.2D& M-10.3C \\
A-2.1D& B-11.2H& A-2.3D& M-10.2C \\
A-2.2D& B-12.1G& A-2.3C& M-10.1C \\
A-2.3D& B-12.3G& A-2.2C& M-10.1B \\
A-2.3C& B-12.2G& A-2.1D& M-10.2B \\
A-2.2C& B-12.3G& A-2.1C& M-10.3B \\
A-2.1C& B-12.1H& A-2.1B& M-10.3A \\
A-2.1B& B-12.2H& A-3.2A& M-10.2A \\
A-2.2B& B-12.3H& A-2.1A& M-10.1A \\
&       B-12.3J& A-2.2B& M-10.2A \\
&       B-12.2J& A-2.3B& P-10.1A \\
&       B-12.1J& A-2.3A& P-10.2A \\
&       B-11.1J& A-3.1A& P-10.3A \\
&       B-11.2J& A-2.2A& P-10.3B \\
&       B-11.3J& A-2.3A& P-10.2B \\
&       B-11.3H&       & P-10.1B \\
&       B-11.3J&       & P-10.1C \\
&       B-11.2G&       & P-10.2C \\
&              &       & P-10.3C \\
&              &       & P-10.3D \\
\hline
\end{tabular}
\end{center}
}
\end{table}

The lists have patterns, suggesting drivers are following higher-level clusters created by transitions in the four components of the IDs.
To exploit this, we create such a clustering by selecting a subset ${\cal S}$ of the symbols $\{\Gamma, x, y, \Delta\}$ and grouping all IDs that have matching values in the ${\cal S}$ positions.
{\em Super clusters} are determined by a selection ${\cal S}$ of three symbols, {\em super-super clusters} are determined by a selection ${\cal T} \subset {\cal S}$ of two symbols, and a {\em top-level clustering} is determined by ${\cal U} \subset {\cal T}$ having a single symbol.
For example, setting ${\cal S} = \{\Gamma, x, \Delta\}$, ${\cal T} = \{\Gamma, x\}$, and ${\cal U} = \{\Gamma\}$ gives for route amz0002 from Table~\ref{zone_table} the super clusters $\{\mbox{A-2.2E, A-2.1E}\}$, $\{\mbox{A-2.1D, A-2.2D, A-2.3D}\}$, $\{\mbox{A-2.3C, A-2.2C, A-2.1C}\}$, $\{\mbox{A-2.1B, A-2.2B}\}$ and a single super-super cluster consisting of all of the zones.
The choices of ${\cal S}$, ${\cal T}$, and ${\cal U}$ are made in a preliminary  computation, minimizing the number of times the training tours cross the components of the partitions of zones at each clustering level.

The super clusters form a partition of the zone clusters, and the super-super clusters form a partition of the super clusters.
We can therefore handle super-cluster constraints and super-super-cluster constraints with penalties, as we discussed in the case of zone clusters.
Computational tests with these new constraints are reported in Table~\ref{table_super}.
\begin{table}[htb]
\caption{Scores of LKH-AMZ runs with multi-level clusters (M1 processor)}
\label{table_super}
\begin{center}
\begin{tabular}{c|ccc} \hline
Run Time & Clustering &  +Precedence & +Transitive \\ \hline
1s       & 0.03817    & 0.02715      & 0.02801 \\
2s       & 0.03815    & 0.02718      & 0.02793 \\
5s       & 0.03820    & 0.02716      & 0.02794 \\ \hline
\end{tabular}
\end{center}
\end{table}
The ``Clustering" column gives results with only the multi-level clusters; the scores  are similar to a 0.03830 score we obtained using optimal tours for these instances provided by Concorde (modeling cluster constraints with transformed travel times).
Much better scores are obtained by adding component-path precedence constraints, reported in the ``+Precedence" column, or the transitive-closure precedence constraints, reported in the ``+Transitive" column.

Going further with the analysis, consider again our example route amz0002.
The zone order for the route is displayed in Table~\ref{table_amz0002}, with one super cluster in each row.
\begin{table}[htb]
\caption{Zone order for route amz0002 arranged by super clusters}
\label{table_amz0002}
\begin{center}
\begin{tabular}{c|ccc} \hline
Super Cluster  & Zone 1   & Zone 2  & Zone 3  \\ \hline
E              &          & A-2.2E  & A-2.1E  \\
D              &  A-2.1D  & A-2.2D  & A-2.3D  \\
C              &  A-2.3C  & A-2.2C  & A-2.1C  \\
B              &  A-2.1B  & A-2.2B  &         \\ \hline
\end{tabular}
\end{center}
\end{table}

Each super cluster is determined by the $\Delta$ symbol, given in the first column of the table.
Notice that the zone IDs for each super cluster appear in either sorted or reverse sorted order.
For example, the super cluster corresponding to $\Delta = $ D is in sorted order (with symbol $y$ increasing from 1 to 2 to 3), while the super cluster corresponding to $\Delta = $ C is in reverse sorted order (with symbol $y$ decreasing from 3 to 2 to 1).
We can force this behavior in the tours produced by LKH-AMZ by sorting by zone ID the clusters in each super cluster, and adding a neighbor constraint for each pair of zones in the sorted order.

Notice that the final zone A-2.1E in super cluster E has $y$ = 1, matching the $y$ = 1 value for the first zone A-2.1D in the following super cluster D.
Similarly, the final zone A-2.3D in super cluster D has $y = 3$, matching the $y$ = 3 in the first zone A-2.3C in the following super cluster C.
For this example, we can attempt to match the super cluster to super cluster transition by adding neighbor constraints for (A-2.1E, A-2.1D) and for (A-2.3D, A-2.3C).
In general, for neighboring pairs of super clusters $(G,H)$ (after sorting, by super-cluster ID, the super clusters within each super-super cluster) we possibly add neighbor constraints for pairs of zones $(g,h)$, with $g$ being either the first or last zone in the sorted order for $G$ and $h$ being either the first or last zone in the sorted order for $H$.
The neighbor constraint $(g,h)$ is added if these zones have matching values for the unique symbol $u$ in $\{\Gamma,x,y,\Delta\} \setminus {\cal S}$.
If there are two such pairs of zones for $(G,H)$, then we create a disjunction for the pair of neighbor constraints.
Corresponding constraints are added also for sorted orders of super clusters within super-super clusters, and for sorted orders of super-super clusters within each top-level cluster.
Adding the entire collection of new constraints gives the results reported in Table~\ref{table_superadd}.
\begin{table}[htb]
\caption{Scores of LKH-AMZ runs with sorted multi-level clusters (M1 processor)}
\label{table_superadd}
\begin{center}
\begin{tabular}{c|ccc} \hline
Run Time & Sorted Clustering &  +Precedence & +Transitive \\ \hline
1s       & 0.02778           & 0.02047      & 0.02176 \\
2s       & 0.02768           & 0.02045      & 0.02185 \\
5s       & 0.02776           & 0.02045      & 0.02184 \\ \hline
\end{tabular}
\end{center}
\end{table}

To complement the constraints created with zone ID patterns, we add also super-cluster path constraints derived from the set of training routes.
Here we follow the idea adopted in our work on zone-level precedence.
For a routing instance $R$, we find a reference route $Q$ having the greatest number of super clusters in common with $R$.
We add a path constraint for pairs of super clusters $(C,D)$, such that $C$ and $D$ appear in both $R$ and $Q$, the clusters $C$ and $D$ are each entered exactly twice in $Q$, and $C$ and $D$ appear consecutively in $Q$.

Computational results with the super-cluster path constraints are reported in the ``Full" column of Table~\ref{table_final}.
The new constraints are added to the ``+Precedence" model from Table~\ref{table_superadd}, improving the score to under 0.02000.
\begin{table}[htb]
\caption{Scores of LKH-AMZ runs with full and alternate constraints (M1 processor)}
\label{table_final}
\begin{center}
\begin{tabular}{c|cc} \hline
Run Time & Full     &  Alternate \\ \hline
1s       & 0.01993  & 0.02108 \\
2s       & 0.01995  & 0.02107 \\
5s       & 0.01997  & 0.02106 \\ \hline
\end{tabular}
\end{center}
\end{table}
The ``Alternate" column in Table~\ref{table_final} starts with the ``+Transitive" model and adds super-cluster precedence constraints for super-cluster pairs $(C,D)$, rather than super-cluster path constraints.
This gives worse results, but the diversity of the tours produced will prove to be useful in a best-of-two selection described in Section~\ref{section_contest}.

\begin{figure}[htb]
    \centering
  \includegraphics[width=\columnwidth]{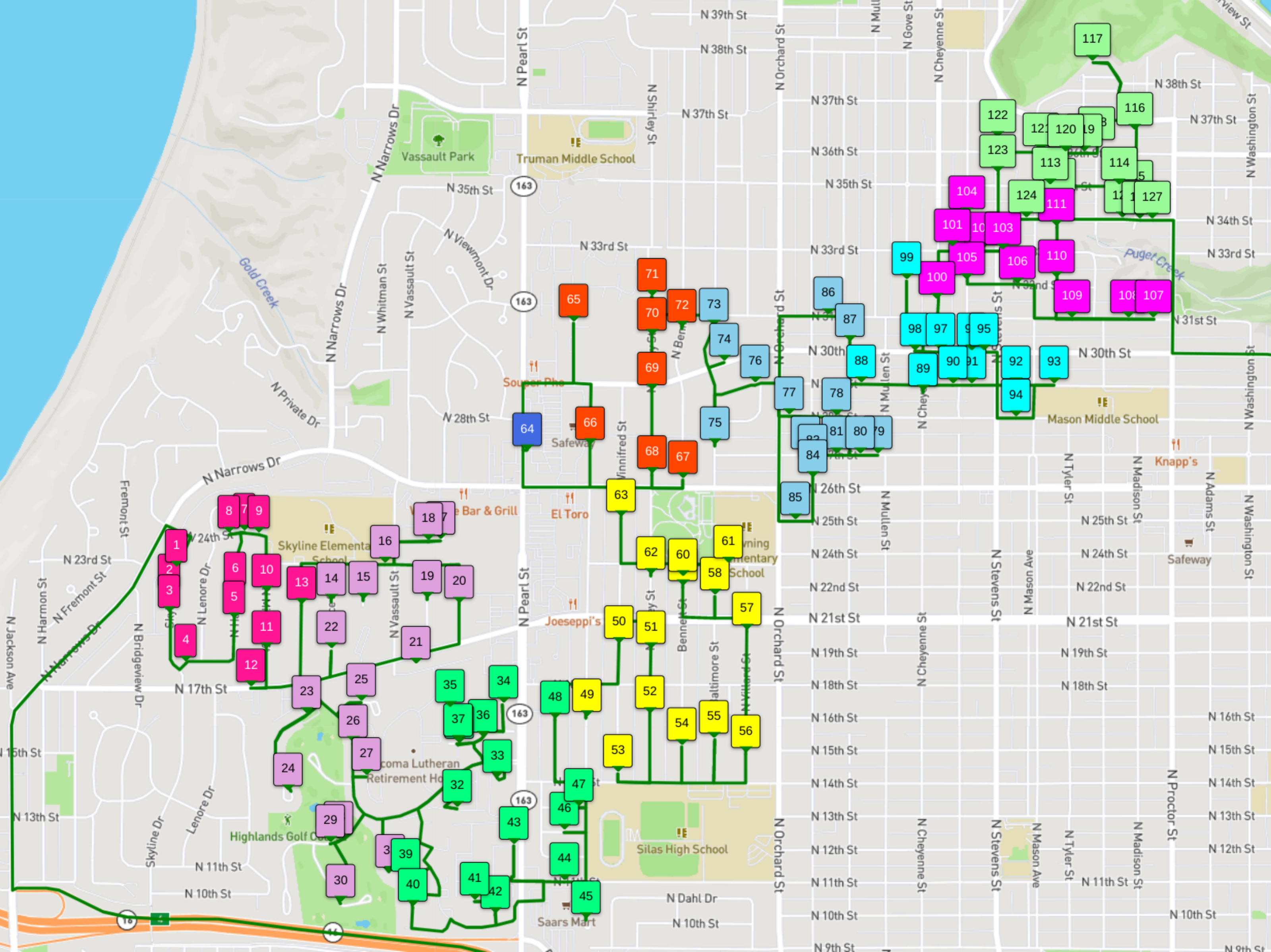}
  \caption{JPT tour for route amz0002. \\~}
  \label{figjpt}
\end{figure}

Figure~\ref{figjpt} shows the solution of  a ``Full'' constraint run on amz0002.
  Note that it serves the zones in the same order as the Amazon driver shown in Figure~\ref{figdriver},
  but it is has different stop orders inside the zones.

As a final remark on time windows, we point out that the ``Full'' constraints with additional time window constraints yield a significantly worse score of  0.02148  in a 1s run.

\subsection{Complexity of a penalty computation}\label{section_complexity_penalties}
To compute the penalty of a tour, we perform a single pass through the tour (an $O(n)$ operation).
During the search for improving moves it would be desirable to have a computation of lower complexity.
We explored the possibility of computing the penalty change  for a provisional move.
Possibly, this could be computed by considering only those stops that would be touched by the move.
However, it was not clear how to implement it and, in our experiments on the benchmark instances, it did not turn out to be a problem.
\section{Competition results}
\label{section_contest}

The Last Mile Challenge evaluations were run on an AWS EC2 m5.4xlarge server, with codes given a time limit of 12 hours for a build phase (where information is extracted from training instances; during this phase codes do not have access to the 3,072 new instances) and 4 hours for an apply phase (where tours are computed).
The AMS server has an 8-core processor and supports 16 virtual cores.

The large amount of computing power assigned to the build phase supported the use of machine-learning techniques, such as training a deep neural network.
In our submission, we used this time only to read and copy the training data, fill in missing zone IDs, select zone ID patterns for multi-level clusters, and build the component paths.
We elected to not fit the penalty values to the training instances, using instead only two levels: zone-precedence constraints and zone-neighbor constraints receive a penalty of 1 when violated, and zone-path constraints and all super-cluster and super-super-cluster constraints receive the penalty 1000.
The total time in the build phase was 109 seconds on a single core of the AWS server.

\subsection{Best-of-two tours}
\label{subsection_merge}

LKH-AMZ is designed to very quickly find good tours satisfying all or most of the specified constraints.
Indeed, as can be seen in the results reported in Table~\ref{table_final}, when looking to minimize the score, it is not always productive to increase the LKH-AMZ time limit.
In the competition, we instead split the computation, running LKH-AMZ on both the full and alternate models (from Table~\ref{table_final}).

While the full model usually gives  better scores than the alternate model (see Table~\ref{table_final}),
  there are  also outliers where  the alternate model gives a better score.
  Unfortunately, it is difficult to predict where the alternate solution has a better score.
  In a best-of-two strategy, we choose the alternate solution if its length is substantially shorter.

More precisely, we compute the travel times $t_{f}$ and $t_{a}$ for an  instance's tours computed with the full and the alternate constraints, then select the full tour if $t_f \leq MergeFactor * t_a$, and otherwise select the alternate.
We use $MergeFactor = 1.01$  to replace length outlier tours only.

The best-of-two results are presented in Table~\ref{table_merge}.
\begin{table}[htb]
\caption{Scores of best-of-two LKH-AMZ tours with full and alternate constraints (M1 processor)}
\label{table_merge}
\begin{center}
\begin{tabular}{c|ccc} \hline
 Run Time        & Full    & Alternate & Best-of-two  \\ \hline
 1s+1s           & 0.01993 & 0.02018   & 0.01978 \\
 2s+1s           & 0.01995 & 0.02018   & 0.01979 \\
 2s+2s           & 0.01995 & 0.02107   & 0.01980 \\
 5s+2s           & 0.01997 & 0.02107   & 0.01980 \\
 5s+5s           & 0.01997 & 0.02106   & 0.01980 \\ \hline
\end{tabular}
\end{center}
\end{table}
The ``Run Time" column indicates the time allocated to individual routes, where the first value is the LKH-AMZ time for the full model and the second value is the time for the alternate model.
The 1s+1s best-of-two gives our best result, producing tours for the 1,107 High+Delivered instances having a mean score of 0.01978 and a median score of 0.00732.

\subsection{Evaluation run}

The version of LKH-AMZ presented in this paper contains a number of improvements over the code submitted to the competition, particularly in terms of computation speed.
Thus, for the competition we elected to use a much longer run time than we considered in Table~\ref{table_merge}.

The {\tt model\_apply} script in our submitted code determines the amount of time to allocate to LKH-AMZ after the constrained ATSP instances are created, aiming for a total run time of 3.5 hours using 16 threads.
In a test on an AWS EC2 m5.4xlarge server using 3,050 routes (selected from the training instances), the code adopted a 41s+20s best-of-two merge.

The top three results in the competition are listed in Table~\ref{table_scores}.
Our submission, {\em Just Passing Through}, received the \$100,000 prize for a score of 0.0248.
\begin{table}[htb]
\caption{Competition scores of top three teams}
\label{table_scores}
\begin{center}
\begin{tabular}{ccc} \hline
Team                   & Score   & Prize     \\ \hline
{\em Just Passing Through}                   & 0.0248  & \$100,000 \\
Xiaotong Guo, Baichuan Mo, Qingyi Wang & 0.0353  & \$50,000  \\
Okan Arslan, Rasit Abay              & 0.0391  & \$25,00   \\ \hline
\end{tabular}
\end{center}
\end{table}
We remark that the second and third placed teams also adopted clustered ATSP models; short reports on their work can be found in the challenge proceedings (\cite{lmrr2021}).

\subsection{Computer implementation}
LKH-AMZ contains approximately 3,500 lines of C code.
It links only the standard C library, making it easily portable to many computing platforms.
The input to LKH-AMZ consists of  a distance matrix and a specification of zones and zone constraints.
  It reads files in the TSPLIB format, which we extended by the constraints explained in Section~\ref{section_constraints}.

Together with the LKH-AMZ code, our submission contains Python scripts, totaling approximately 2,000 lines, used for analyzing zone information, extracting constraints, and visualizing tours.
An additional 500 lines of C implement an internal scoring function, 300 lines of shell scripts control the LKH-AMZ solver and execute the best-of-two routine, and 100 lines of scripts control the build and apply phases on the AWS server.
All codes and scripts are available under the MIT License at {{\url{https://github.com/heldstephan/jpt-amz}}}.

\subsection{Reduced variants of LKH-AMZ}
We made runs with several variants of LKH-AMZ, omitting certain features to assess their relative importance.
The results are shown in Table~\ref{tab:reduced-effort}.
The computations were made on a MacBook Pro equipped with an M1 Pro processor and the time limit was set to 1 second per instance.

The first line in the table shows results for the standard version of LKH-AMZ described in this paper, adopting 3-opt and 4-opt moves, 8\textit{n} trials, and as many runs as the time limit allows.
In the second line we report results obtained with 4-opt disabled, using 3-opt moves only.
This raised the score only slightly, from 0.01981 to 0.02011.
In the subsequent experiments we omitted 4-opt moves.

Next, we turned off the  iterated local search from Section~\ref{sec_iterated_ls} by restricting MAX\_TRIALS to 1, while still performing multiple  runs from new random starting solutions (see Section~\ref{sec_multiple_runs}).
This raised the score further to 0.02243.
Going the other way, using 8\textit{n} trials but only a single global run gives a worse score of 0.02724.
Turning off both options, that is, a single trial using 3-opt moves, gives a very poor score of 0.07135.

\begin{table}[htb]

  \begin{center}
    \caption{Results of different variants omitting 4-opt moves, kicks, or multiple runs in LKH-AMZ using a 1s time limit (M1 Pro processor).}
    \label{tab:reduced-effort}

      \begin{tabular}{cccc}
      MOVE\_TYPE & MAX\_TRIALS &  Multiple Runs   & Score\\\hline
      3-opt \& 4-opt &  8\textit{n} &   yes   & 0.01981\\\hline
      3-opt          &  8\textit{n} &  yes    & 0.02011\\
      3-opt          &        1      &    yes  & 0.02243\\
      3-opt          &  8\textit{n} &      no     & 0.02724\\
      3-opt          &      1        &     no      & 0.07135\\
      \end{tabular}

  \end{center}

\end{table}

\subsection{Real-time optimization on servers and handhelds}
\label{section_practice}

The previous experiments have demonstrated that the code is fast enough to  (re-)compute tours
ad-hoc on a driver's handheld device within a second.
Table~\ref{table_short} shows that we could even go into the sub-second range for ``real-time'' computations.
Already with a  0.1 second  running time limit the code would presumably have won the Amazon contest.
However, there is a big jump in the quality when increasing the running time to 0.2 seconds.
\begin{table}[htb]
\caption{Short LKH-AMZ runs on 1,107  instances with full constraints (M1 processor)}
\label{table_short}
\begin{center}
\begin{tabular}{c|ccc} \hline
Run Time & Score   \\ \hline
0.1s     & 0.02488 \\
0.2s     & 0.02042 \\
0.3s     & 0.02007 \\
0.4s     & 0.02000 \\
0.5s     & 0.01995 \\ \hline
\end{tabular}
\end{center}
\end{table}

In settings where it is necessary to process a large number of routes in a short amount of time, LKH-AMZ can be run in batch mode.
Here, many jobs compete for the limited memory bandwidth of a server, in particular when enforcing hyper-threading.
  A job could be idling while its CPU clock is running, degrading the quality.
Using the included {\tt solve} script, the number of parallel threads can be set, together with the running time per instance.
In Table~\ref{table_parallel} we display batch results using full constraints without best-of-two selection for faster running times.
We conducted all runs on a Linux server equipped with two 16-core Intel Xeon Gold 5218 CPU @ 2.30GHz processors and on a Linux server equipped with two 64-core AMD EPYC2 7742 processors.
We once assigned as many parallel jobs as there are processing cores, and once  we used twice that number. The second mode assigns two jobs  on each  processor, i.e. it is enforcing  hyper-threading. The wall clock times show the total time for solving all 1,107 instances. The one second  runs on the AMD server are dominated by reading and writing the instances, which is not considered in the
  running time limit.  It can be seen that hyper-threading does not degrade the results.
The AMD server, in particular, gives impressive results, delivering tours with scores under 0.02 at a rate of over 6,000 instances per minute of wall clock time.
\begin{table}[htb]
\caption{LKH-AMZ batch processing on 1,107 instances with full constraints}
\label{table_parallel}
\begin{center}
\begin{tabular}{c|ccc|ccc} \hline
\multicolumn{1}{c}{ } & \multicolumn{3}{c}{Intel Xeon Gold} & \multicolumn{3}{c}{AMD EPYC2 7742}\\ \hline
Run Time & Threads &  Score   & Wall Clock & Threads & Score & Wall Clock \\ \hline
1s       &  64     & 0.02010  &  20s  & 256 & 0.01997 & 10s \\
1s       &  32     & 0.02003  &  37s  & 128 & 0.01997 & 11s \\
2s       &  64     & 0.01993  &  38s  & 256 & 0.01993 & 13s \\
2s       &  32     & 0.01993  &  73s  & 128 & 0.01993 & 20s \\ \hline
\end{tabular}
\end{center}
\end{table}

\section{Conclusions\label{section_conclusions}}

We presented how to efficiently model and solve complex constraints in the ATSP.
These constraints have many applications beyond the Amazon competition: precedence constraints occur in pickup and delivery problems, or as disjunctive  constraints in the cable tree wiring problem.
Hierarchical clusters occur whenever parcels are stored in containers, carts, bags, etc.

Furthermore, we described new learning and constraint extraction methods to reflect information from previously driven routes in newly planned routes.

The LKH-AMZ algorithm computes high-quality solutions for these constrained instances within less than a second.
  Thus, it meets the demands of computing thousands of tours in the planning period of a real-world logistics application.
It is designed to handle stop-level constraints, such as requiring stop $a$ to precede stop $b$ in a tour.
Analysis of larger data sets (with repeated visits to the same streets) could make use of this option, capturing local preferences of van drivers.

Concerning the Amazon benchmark data, it remains difficult to assess to which extent our approach recovers
the driver behavior, possibly deviating from the planned tour, or just the Amazon planning process.
To measure this, it would be necessary to know the undisclosed tours that were proposed to the drivers.

\bigskip

{\RaggedRight
\bibliographystyle{plainnat}
\bibliography{jpt}
}

\end{document}